\documentclass[12pt]{amsart}
\usepackage{amsfonts}
\usepackage{amssymb}
\theoremstyle{definition}

\theoremstyle{remark}

\begin{document}

\centerline{\large\bf ON THE AXIOM OF PLANES AND THE AXIOM}
\centerline{\large\bf OF SPHERES IN THE ALMOST HERMITIAN GEOMETRY
\footnote{\it SERDICA Bulgaricae mathematicae publicationes. Vol. 8, 1982, p. 109-114}}

\vspace{0.2in}
\centerline{\large OGNIAN T. KASSABOV}

\vspace{0.5in}
{\sl We prove analogues of Cartan's criterion for constancy of sectional curvature
to an arbitrary almost Hermitian manifold. As a consequence we establish for such a 
manifold analogues of a Cartan's theorem. Our results generalize some theorems
in [2; 5; 7, 8].}

\vspace{0.2in}
{\bf 1. Introduction.} Let $N$ be an $n$ - dimensional submanifold of an $m$-dimensional
Riemannian manifold $M$ with Riemannian metric $g$ and let $\tilde\nabla$ and $\nabla$ 
be the Levi-Civita connections on $M$ and $N$, respectively. It is well known, that 
the equation 
$$
	\alpha(x,y) = \widetilde\nabla_xy-\nabla_xy \ ,
$$
where $x,\ y \in {\mathfrak X}N$, defines a normal-bundle-valued symmetric tensor field, called
the second fundamental form of the immertion. The submanifold $N$ is said to be totally
umbilical, if $\alpha(x,y)=g(x,y)H$ for all $x,\ y \in {\mathfrak X}N$, where $N=(1/n){\rm trace}\, \alpha$
\ is the mean curvature vector of $N$ in $M$. In particular, if $\alpha$ vanishes identically,
$N$ is called a totally geodesic submanifold of $M$.

For $x \in \mathfrak XN,\ \xi \in \mathfrak XN^{\perp}$ we write 
$ \widetilde\nabla_x\xi =- A_{\xi}x+D_x{\xi}$, where $-A_{\xi}x$ (respectively, $D_x\xi$)
denotes the tangential (respectively, the normal) component of $\widetilde\nabla_x\xi$. A normal
vector field $\xi$ is said to be parallel, if \ $D_x\xi=0$ \ for each \ $x\in \mathfrak XN$.

The manifold $M$ is said to satisfy the axiom of $n$-planes (respectively, $n$-spheres),
where $n$ is a fixed integer $2 \le n<m $ if for each point $p\in M$ and for any 
$n$-dimensional subspace $\alpha$ of $T_pM$ there exists an n-dimensional totally
geodesic submanifold $N$ (respectively an $n$-dimensional totally umbilical submanifold
$N$ with non-zero parallel mean curvature vector) containing $p$, such that
$T_pN=\alpha$.

In his book on Riemannian geometry [1] E. C\,a\,r\,t\,a\,n proved the following theorem.

\vspace{0.05in}
T\,h\,e\,o\,r\,e\,m. {\it Let $M$ be an $m$-dimensional Reimannian manifold, $m>2$,
which satisfies the axiom of planes for some $n$, $2\le n < m$. Then $M$ has
constant sectional curvature. }

\vspace{0.05in}
In [4] L\,e\,u\,n\,g and N\,o\,m\,i\,z\,u have substituted the axiom of $n$-planes with
the axiom of $n$-spheres and have proved a generalization of the above mentioned
Cartan's theorem.

Analogous results for Kaehler manifolds have been proved in [2;\ 5;\ 8] and in
[7] it has been studied a similar problem for some almost Hermitian manifolds.

\vspace{0.2 in}
{\bf 2. Preliminaries.} Let $M$ be an $m$-dimensional Riemannian manifold with
Riemannian metric $g$ and let $\nabla$ be its Levi-Civita connection. The 
curvature tensor $R$ associated with $\nabla$ has the following properties:

1) $ R(X,Y)=-R(Y,X)$ \ for $X,\ Y \in T_pM$,

2) $ R(X,Y)Z+R(Y,Z)X+R(Z,X)Y=0$ \ for $X,\ Y,\ Z \in T_pM$,

3) $ R(X,Y,Z,U)=-R(X,Y,U,Z)$ \ for $X,\ Y,\ Z,\ U \in T_pM$,

\noindent where $ R(X,Y,Z,U)=g(R(X,Y)Z,U)$. 
 
The curvature of a two dimensional plane in $T_pM$ with an orthonormal
basis $X,\ Y$ is defined by $K(X,Y)=R(X,Y,Y,X)$.

It is easy to compute, that if \ $N$ \ is a totally geodesic submanifold of \ $M$ \
or a totally umbilical submanifold of \ $M$ \ with parallel mean curvature
vector, then $R(x,y,z,\xi)=0$ \ for all vectors \ $x,\, y,\, z \in T_pM$,
$\xi \perp T_pM$ \ and for each point \ $p\in N$.

Now, let \ $M$ \ be a $2m$-dimensional almost Hermitian manifold with Riemannian
metric \ $g$ \ and almost complex structure \ $J$. 

A subspace  $\alpha$  in  $T_pM$  is said to be holomorphic (respectively,
antiholomorphic or totally real) if  $J\alpha = \alpha$  (respectively  
$J\alpha \perp \alpha)$. For the dimension \ $k $ \ of a  holomorphic
(respectively, antiholomorphic) subspace  $\alpha$  of  $T_pM$  we have  
$k=2n$,  $1\le n\le m$  (respectivey, $1\le k \le m$). If the holomorphic
(respectively, antiholomorphic) sectional curvature in each point $p\in M$,
i.e. the curvature of a holomorphic (respectively, antiholomorphic) 2-dimensional subspace 
$\alpha$  of  $T_pM$  does not depend on  $\alpha$, then  $M$  is said
to be of pointwise constant holomorphic (respectively, antiholomorphic)
sectional curvature. 

A connected Riemannian (respectively, Kaehler) manifold of global constant 
sectio\-nal curvature (respectively, of constant holomorphic sectional curvature)
is called a real-space-form (respectively, a complex-space-form).

An almost Hermitian manifold is said to be an $RK$-manifold, if \ 
$R(X,Y,Z,U)=R(JX,JY,JZ,JU)$ \  for all \ $X,\, Y,\, Z,\, U \in T_pM$, \ $p\in M$. 

For a two dimensional subspace \ $\alpha$ \ of \ $T_pM$ \ with an orthonormal basis
$X,\ Y$ the angle \ $ \theta \in [0,\pi/2]$ \ between \ $\alpha$ \ and \ $J\alpha$ \ 
is defined by \ $\cos\theta = |g(X,JY)|$.

We shall need the following theorems: 

\vspace{0.05in}
T\,h\,e\,o\,r\,e\,m \, A [3]. {\it Let \ $M$ \ be a \ $2m$-dimensional almost Hermitian manifold,
$m\ge 2$, and let \ $T\, : \, (T_pM)^4 \rightarrow \bf R$ \ be a four-linear mapping, 
which satisfies the conditions:

1) for all \ $X,\, Y,\, Z,\ U \in T_pM$

$
	\qquad\qquad	T(X,Y,Z,U)=-T(Y,X,Z,U) \ ,  
$

$
	\qquad\qquad	T(X,Y,Z,U)+T(Y,Z,X,U)+T(Z,X,Y,U)=0 \ ,
$

$
	\qquad\qquad	T(X,Y,Z,U)=-T(X,Y,U,Z) \ ; 
$ 

2) $T(X,Y,Y,X)=0$ , where $X,\ Y$ is a basis of an arbitrary two dimensional subspace \
$\alpha$ \ in \ $T_pM$, for which the angle between \ $\alpha$ \ and \ $J\alpha$ \ is
one of the numbers \ $0, \ \pi/4,\ \pi/2$.

Then $T=0$.} 

\vspace{0.05in}
Let for all \ $X,\, Y,\, Z,\ U \in T_pM$
$$
	R_1(X,Y,Z,U)=g(X,U)g(Y,Z)-g(X,Z)g(Y,U) \ , \qquad\qquad\qquad\qquad
$$
$$
	R_2(X,Y,Z,U)=g(X,JU)g(Y,JZ)-g(X,JZ)g(Y,JU)-2g(X,JY)g(Z,JU) \ .
$$

\vspace{0.05in}
T\,h\,e\,o\,r\,e\,m \ B [3]. {\it If \ $M$ \ is a \ $2m$-dimensional \ $RK$-manifold, $m\ge 2$,
with pointwise constant holomorphic sectional curvature \ $c$ \ and with pointwise 
constant antiholomorphic sectional curvature \ $K$, then the curvature tensor \ 
$R$ \ has the form} 
$$
	R=KR_1+(c-K)R_2/3 \ .  \leqno (2.1)
$$

\vspace{0.05in}
As is proved in [6], if the curvature tensor of a $2m$-dimensional connected
almost Hermitian manifold has the form (2.1) and if  $m \ge 3$, then  $c $  and  $K$  
are global constants. On the other hand, it is proved in [3], that if the curvature 
tensor of an almost Hermitian manifold  $M$  of dimension  $2m \ge 4$  has the
form (2.1) with global constants  $c$  and  $K$, then either  $M$  is of 
constant holomorphic sectional curvature  $c=K$  or  $M$ is a Kaehler manifold of 
constant holomorphic sectional curvature. Hence we have:

\vspace{0.05in}
T\,h\,e\,o\,r\,e\,m \ C. {\it Let \ $M$ \ be a connected \ $RK$-manifold of dimension \ 
$2m \ge 6$. If \ $M$ \ has pointwise constant holomorphic sectional curvature and 
pointwise constant antiholomorphic sectional curvature, then \ $M$ \ is one of the 
following:

1) a real-space-form;

2) a complex-space-form.}   

\vspace{0.2in}

{\bf 3. Criterions for constancy of the holomorphic and the antiholomorphic curvature
at one point.}

\vspace{0.05in}
L\,e\,m\,m\,a \ 1. {\it Let \ $M$ \ be an almost Hermitian manifold with dimension \ $2m$,
$m\ge 2$ \ and for a point \ $p\in M$
$$
	R(X,JX,JX,Y)=0          \leqno (3.1)
$$
holds for all \ $X,\, Y \in T_pM$, with \ $g(X,Y)=g(X,JY)=0$. Then \ $M$ \ has constant
holomorphic sectional curvature at \ $p$ \ and
$$
	R(X,Y,Y,X)=R(JX,JY,JY,JX)  \ ,     \leqno (3.2)
$$
where $X,\ Y$ are as above.}

P\,r\,o\,o\,f. Taking two arbitrary unit vectors  $X, Y$  in  $T_pM$  with 
  $g(X,Y)=g(X,JY)=0$  and applying (3.1) for the vectors \ $X+\alpha Y$,
$\alpha X-Y$, we obtain 
$$
	\begin{array}{c}
		H(X)-\alpha^2H(Y)+(\alpha^2-1)R(X,JX,JY,Y)+(\alpha^2-1)R(X,JY,JX,Y) \\
		+\alpha^2K(X,JY)-K(JX,Y)=0 \ ,  
	\end{array}                   \leqno (3.3)
$$
where  \ $H(X)=R(X,JX,JX,X)$ \ denotes the holomorphic sectional curvature,
determined by \ $X$.

Let \ $\alpha = 1$:
$$
	H(X)-H(Y)+K(X,JY)-K(JX,Y)=0  \ .           \leqno (3.4)
$$

From (3.3) and (3.4) it follows
$$
	H(Y)=R(X,JX,JY,Y)+R(X,JY,JX,Y)+K(X,JY)   \ .    \leqno (3.5)  
$$

Analogously \ $H(X)=R(X,JX,JY,Y)+R(X,JY,JX,Y)+K(JX,Y)$.

Substituting \ $X$ \ by \ $JX$ \ and \ $Y$ \ by \ $JY$ we  get
$$
	H(X)=R(X,JX,JY,Y)+R(X,JY,JX,Y)+K(X,JY) \ .    \leqno (3.6)
$$

From (3.5) and (3.6) we see that 
$$
	H(X)=H(Y)      \leqno (3.7)
$$
and combining this with (3.4) we find (3.2).

Let \ $m>2$ \ and \ $U,\, V$ \ be arbitrary unit vectors in \ $T_pM$.
We choose \ $X$ in ${\rm span} \{ U,JU \}^{\perp} \cap {\rm span} \{ V,JV \}^{\perp}$.
According to (3.7) we have \ $H(U)=H(X)=H(V)$ \ and the Lemma is proved in the case
\ $m>2$. In the case \ $m=2$ \ we put \ $c=H(X)=H(Y)$ \ and using (3.6) we see that
\ $H(\alpha X+\beta Y)=c$, where $\alpha^2+\beta^2=1$. Hence it is not difficult to 
find that the holomorphic sectional curvature in \ $p$ \ is a constant.

\vspace{0.05in}
The following lemma is trivial.

\vspace{0.05in}
L\,e\,m\,m\,a \ 2. {\it Let \ $\alpha$ \ be a two dimensional subspace in \ $T_pM$ \ such
that the angle between \ $\alpha$ \ and \ $J\alpha$ \ is \ $\pi/4$. Then \ $\alpha$ \
has an orthonormal basis \ $X, \, (JX+U)/\sqrt 2$, where \ $X,\, U$ \ are unit vectors 
in \ $T_pM$ \ with \ $g(X,U)=g(X,JU)=0$.}  

\vspace{0.05in}
L\,e\,m\,m\,a \ 3. {\it Let \ $M$ \ be a \ $2m$-dimensional almost Hermitian manifold, \ 
$m\ge2$ \ and for each point \ $p\in M$ \ (3.1) holds for all \ $X,\, Y \in T_pM$ \
with \ $g(X,Y)=g(X,JY)=0$. Then \ $M$ \ is an \ $RK$-manifold.}

P\,r\,o\,o\,f. We put \ $T(X,Y,Z,U)=R(X,Y,Z,U)-R(JX,JY,JZ,JU)$ \ for all \ 
$X,\, Y,$ $ Z,\, U \in T_pM$. Obviously \ $T$ \ has the property 1) of Theorem A.
Let \ $\alpha$ \ be a subspace in \ $T_pM$ \ such that the angle between \ 
$\alpha$ \ and \ $J\alpha$ \ is \ $\theta$. If \ $\theta = 0$, $\alpha$ \ is a
holomorphic plane and if \ $\alpha= {\rm span}\{ X,JX \}$ \ we have \ 
$ T(X,JX,JX,X)=0$. If \ $\alpha = \pi/2$, $\alpha$ is an antiholomorphic plane 
and we can choose two vectors \ $X,\, Y \in T_pM$ \ such that 
$\alpha= {\rm span}\{ X,JY \}$,  $g(X,Y)=g(X,JY)=0$. According to Lemma 1
we have \ $T(X,Y,Y,X)=0$. Let \ $\theta=\pi/4$ \ and let  \ $X, \, (JX+U)/\sqrt 2$
be an orthonormal basis of \ $\alpha$ \ as in Lemma 2. Then \
$T(X,JX+U,JX+U,X)=0$. According to Theorem A we have \ $T=0$, which proves our
assertion.

\vspace{0.05in}
L\,e\,m\,m\,a \ 4. {\it Let \ $M$ \ be a \ $2m$-dimensional almost Hermitian manifold, \ 
$m\ge2$ \ and for a point \ $p\in M$
$$
	R(X,Y,Y,Z)=0                 \leqno (3.8)
$$ 
holds for all \ $X,\, Y,\, Z \in T_pM$ \ with \ $g(X,Y)=g(X,JY)=g(X,Z)=g(Y,Z)=0$.
Then \ $M$ \ has constant antiholomorphic sectional curvature at \ $p$.}

P\,r\,o\,o\,f. According to Lemma 1 \ $M$ has constant holomorphic sectional curvature \ 
$c$ \ at \ $p$. We apply (3.8) for the vectors \ $X+JX,\, Y, \ JX-X$, where    \
$X,\, Y \in T_pM$ are arbitrary unit vectors with  \ $g(X,Y)=g(X,JY)=0$ and we get \
$K(JX,Y)=K(X,Y)$. As in the proof of Lemma 1 we have
$$
	H(X)=R(X,JX,JY,Y)+R(X,JY,JX,Y)+K(X,JY) \ .    \leqno (3.9)
$$ 

Hence applying the first Bianchi's identity we obtain
$$
	H(X)=2R(X,JX,JY,Y)+R(JX,JY,X,Y)+K(X,JY) \ .    \leqno (3.10)
$$

The substitution of \ $Y$ \ with \ $JY$ \ in (3.9) gives
$$
	H(X)=R(X,JX,JY,Y)-R(X,Y,JX,JY)+K(X,Y) \ .    
$$

Combining this with (3.10) we derive
$$
	2H(X)=3R(X,JX,JY,Y)+K(X,Y)+K(X,JY) \ .    \leqno (3.11)
$$

Let \ $m=2$. We put \ $K=K(X,Y) $ \ and from (3.9), (3.10), (3.11) we have
$ R(X,JX,JY,Y)=\frac 23\{ c-K \}$; $R(X,JY,JX,Y)=\frac 13\{ c-K \}$;
$ R(JX,JY,X,Y)=\frac 13 \{K-c\}$.

We put \ $R'=KR_1+\frac{c-K}3 R_2$. A simple calculation shows that
$ R(X_1,X_2,X_3,X_4)=R'(X_1,X_2,X_3,X_4)$, whenever $X_1,\, X_2,\, X_3,\, X_4$ \ 
are choosen among the vectors \ $X,\, JX,$ $ JY,\, Y$. Consequently \ $R=R'$ \ 
and the Lemma is proved in the case \ $m=2$.

Let \ $m>2$. We choose a unit vector \ $Z$, normal to \ $X,\ JX,\, Y,\, JY$. 
Because of (3.8), from \ $R(X+Z,Y,Y,X-Z)=0$ \ we get
$$
	K(X,Y)=K(Y,Z)    \ .    \leqno (3.12)
$$

Let \ $m=3$. We shall show that
$$
	R(X,JX,Y,Z)=R(X,Y,Z,JX)=0      \leqno (3.13)
$$
and the case \ $m=3$ \ will follow as the case \ $m=2$. From \ 
$ R(\alpha X+JZ,\alpha JX-Z,\alpha JX-Z,Y)=0$, where \ $\alpha$ \ 
takes the values 1 and $-1$, we find
$$
	R(X,JX,Z,Y)+R(X,Z,JX,Y)=0      \leqno (3.14)
$$
and from \ $ R(Y,X+JX,X+ JX,Z)=0$ it follows
$$
	R(X,Y,Z,JX)+R(X,Z,Y,JX)=0      \ .         \leqno (3.15)
$$ 

Using (3.14), (3.15) and the properties of the curvature tensor we get \
$R(X,Y,Z,JX)=0$ \ and together with (3.14) this gives (3.13).

Now let \ $m>3$. We take arbitrary antiholomorhic spaces \ $\alpha,\, \beta$ \
in \ $T_pM$ \ with orthonormal basis \ $X,\, Y$ \ and \ $Z,\, U$ \ respectively 
such that \ $X \perp Y,\, JY$ \ and \ $Z \perp U,\, JU$. Let \ $V,\, W$ \ 
be unit vectors in ${\rm span}\{ X,JX \}^{\perp} \cap  {\rm span}\{ Z,JZ \}^{\perp}$
\ and \ $V \perp W,\, JW$. According to (3.12)
$$
	K(X,V)=K(V,W)=K(V,Z)  \ .     \leqno (3.16)
$$
Let \ $A\in {\rm span}\{ V,JV \}^{\perp} \cap  {\rm span}\{ Z,JZ \}^{\perp} \cap  {\rm span}\{ U,JU \}^{\perp}$
\ be a unit vector. From (3.12)
$$
	K(V,Z)=K(Z,A)=K(Z,U)  \ .      \leqno (3.17)
$$

Analogously
$$
	K(X,V)=K(X,Y)  \ .      \leqno (3.18)
$$

From (3.16), (3.17), (3.18) it follows \ $K(X,Y)=K(Z,U)$ and the Lemma is proved.

\vspace{0.2in}
{\bf 4. The main results.}
Let \ $M$ \ be a \ $2m$-dimensional almost Hermitian manifold, $m \ge 0$.

\vspace{0.05in}
A\,x\,i\,o\,m  o\,f  h\,o\,l\,o\,m\,o\,r\,p\,h\,i\,c  $2n$-p\,l\,a\,n\,e\,s (respectively,
$2n$-s\,p\,h\,e\,r\,e\,s). {\it For each point $p \in M$ and for any $2n$-dimensional
holomorphic subspace \ $\alpha$ \ of \ $T_pM$ \ there exists a totally geodesic
submanifold \ $N$ \ (respectively, a totally umbilical submanifold \ $N$ \ with
nonzero parallel mean curvature vector) containing \ $p$ \ such that \ $T_pN=\alpha$.} 

\vspace{0.05in}
A\,x\,i\,o\,m \, o\,f \, a\,n\,t\,i\,h\,o\,l\,o\,m\,o\,r\,p\,h\,i\,c  $n$-p\,l\,a\,n\,e\,s (respectively,
$n$-s\,p\,h\,e\,r\,e\,s). {\it For each point $p \in M$ and for any $n$-dimensional
antiholomorphic subspace \ $\alpha$ \ of \ $T_pM$ \ there exists a totally geodesic
submanifold \ $N$ \ (respectively, a totally umbilical sub\-manifold \ $N$ \ with
nonzero parallel mean curvature vector) containing \ $p$ \ such that \ $T_pN=\alpha$.} 

\vspace{0.05in}
T\,h\,e\,o\,r\,e\,m 1. {\it Let \ $M$ \ be a \ $2m$-dimensional almost Hermitian manifold,
$m \ge 2$. If \ $M$ \ satisfies the axiom of holomorphic \ $2n$-planes or the axiom of
holomorphic \ $2n$-spheres for some \ $n$, $1\le n <m$, then \ $M$ \ is an \ 
$RK$-manifold with pointwise constant holomorphic sectional curvature.}

P\,r\,o\,o\,f. The condition gives \ $R(X,JX,JX,Y)=0 $ \ for all vectors \ 
$X,\, Y \in T_pM$ \ with \ $g(X,Y)=g(X,JY)=0$ \ and for each point \ $p \in M$ \ 
and the Theorem follows from Lemma 1 and Lemma 3.

\vspace{0.05in}
T\,h\,e\,o\,r\,e\,m 2. {\it Let  $M$  be a  $2m$-dimensional almost Hermitian manifold,
$m \ge 2$. If  $M$  satisfies the axiom of antiholomorphic \ $n$-planes or the axiom of
antiholomorphic  $n$-spheres for some  $n$, $2\le n \le m$, then  $M$  is an  
$RK$-manifold with pointwise constant holomorphic sectional curvature and with pointwise
constant antiholomorphic sectional curvature and consequently the curvature tensor
has the form (2.1).}

P\,r\,o\,o\,f. By the condition it follows \ $R(X,Y,Y,Z)=0 $ \  
for each point \ $p \in M$ \ and for all \ $X,\, Y,\, Z \in T_pM$ \ with \ 
$g(X,Z)=g(Y,Z)=(X,Y)=g(X,JY)=0$. Now the Theorem follows from Lemmas 1, 3 and 4.

\vspace{0.05in}
By Theorem C and Theorem 2 we derive

\vspace{0.05in}
T\,h\,e\,o\,r\,e\,m 3.{\it Let  $M$  be a  $2m$-dimensional connected almost Hermitian manifold,
$m \ge 3$. If \ $M$ \ satisfies the axiom of antiholomorphic \ $n$-planes or the axiom of
antiholomorphic \ $n$-spheres for some \ $n$, $2\le n \le m$, then \ $M$ \ is one 
of the following:

1) a real-space-form;

2) a complex-space-form.}

\vspace{0.5in}
\centerline{\large R E F E R E N C E S}

\vspace{0.2in}
\noindent
1. E. C\,a\,r\,t\,a\,n. Le\c cons sur la g\'eometrie des espaces de Riemann. Paris, 1946.

\noindent
2. B. -Y. C\,h\,e\,n, K. O\,g\,i\,u\,e. Some characterizations of complex space forms. {\it Duke

 Math. J.}, {\bf 40}, 1973, 797-799.

\noindent
3. G. T. G\,a\,n\,\v{c}\,e\,v. Almost Hermitian manifolds similar to the complex space forms.

{\it C. R. Acad. Bulg. Sci.}, {\bf 32}, 1979, 1179-1182.

\noindent
4. D. S. L\,e\,u\,n\,g, K. N\,o\,m\,i\,z\,u. The axiom of spheres in Riemannian geometry.

{\it J. Different. Geom.}, {\bf 5}, 1971, 487-489.

\noindent
5. K. N\,o\,m\,i\,z\,u. Conditions for constancy of the holomorphic sectional curvature.

{\it J. Different. Geom.}, {\bf 8}, 1973, 335-339.

\noindent
6. F\,r. T\,r\,i\,c\,e\,r\,r\,i, L. V\,a\,n\,h\,e\,c\,k\,e. Curvature tensors on almost Hermitian manifolds.

{\it Trans. Amer. Math. Soc.}, {\bf 267}, 1981, 365-398. 

\noindent
7. L. V\,a\,n\,h\,e\,c\,k\,e. Almost Hermitian manifolds with $J$-invariant Riemann 
curvature

tensor.
{\it Rend. Sem. Math. Univers. Pilitech. Torino}, {\bf 34}, 1975-76, 487-498.

\noindent
8. K. Y\,a\,n\,o, I. M\,o\,g\,i. On real representation of Kaehler manifolds. {\it Ann. Math.}, {\bf 61,}

 1955, 170-189.

\vspace {0.4cm}
\noindent
{\it Center for mathematics and mechanics \ \ \ \ \ \ \ \ \ \ \ \ \ \ \ \ \ \ \ \ \ \ \ \ \ \ \ \ \ \ \ \ \ \
Received 16.7.1980

\noindent
1090 Sofia   \ \ \ \ \ \ \ \ \ \ \ \ \ \ \ \ \  P. O. Box 373}

\end{document}